\documentclass[preprint,12pt]{elsarticle}
\usepackage{mathrsfs}
\usepackage{amsfonts}

\usepackage{graphicx,color}

\usepackage{amssymb}
\usepackage{amsthm}

\usepackage{amsmath,mathrsfs,bm,subfigure,verbatim}

\newtheorem{thm}{Theorem}
\newtheorem{lem}[thm]{Lemma}
\newtheorem{prop}[thm]{Proposition}

\newdefinition{defn}{Definition}
\newdefinition{rmk}{Remark}
\newdefinition{alg}{Algorithm}
\newdefinition{exmp}{Example}
\newproof{pf}{Proof}
\newdefinition{Problem}{Problem}
\newproof{Pf}{Proof of Theorem \ref{maintheo}}
\begin{document}

\begin{frontmatter}




\title{On a C. de Boor's Conjecture in a Particular Case and Related Perturbation}


\author{Zhe Li}
\author{Shugong Zhang\corref{cor1}}
\cortext[cor1]{Corresponding author. +86-431-85165801, Fax
+86-431-85165801} \ead{sgzh@jlu.edu.cn}
\author{Tian Dong}
\address{School of Mathemathics,
Key Lab. of Symbolic Computation and Knowledge Engineering
\textup{(}Ministry of Education\textup{)}, Jilin University,
Changchun 130012, PR China}

\begin{abstract}

In this paper, we focus on two classes of $D$-invariant polynomial subspaces.
The first is a classical type, while
the second is a new class.
With matrix computation, we prove that every
ideal projector with each $D$-invariant subspace belonging to either the first class or the second
is the pointwise limit of Lagrange
projectors.  This verifies a particular case of a C. de Boor's  conjecture
asserting that every complex ideal projector is the pointwise limit
of Lagrange projectors. Specifically, we
provide the concrete perturbation
procedure for ideal projectors of this type.
\end{abstract}


\end{frontmatter}


\section{Introduction}

Polynomial interpolation is to construct a polynomial $p$ belonging
to a finite-dimensional subspace of $\mathbb{F}[\bm{x}]$ from a set
of data that agrees with a given function $f$ at the data set, where
$\mathbb{F}[\bm{x}]:=\mathbb{F}[x_1,\ldots,x_d]$ denotes the
polynomial ring in $d$ variables over the field $\mathbb{F}$. It's important to make the comment that $\mathbb{F}$ is the
real field  $\mathbb{R}$ or the complex field  $\mathbb{C}$ in this paper.
Univariate polynomial interpolation has a well developed theory,
while the multivariate one is very problematic since a multivariate
interpolation polynomial is determined not only by the cardinal but
also by the geometry of the data set, cf. \cite{dBo94,GS00:2}.

Recently,  more and more people are getting interested in ideal interpolation,
which is defined by an \emph{ideal projector} on $\mathbb{F}[\bm{x}]$,
namely a linear idempotent operator on $\mathbb{F}[\bm{x}]$ whose
kernel is an ideal, cf. \cite{Bir1979}.
When the kernel of an ideal projector $P$ is the vanishing ideal of
certain finite set $\Xi$ in $\mathbb{F}^d$, $P$ is a \emph{Lagrange
projector} which
provides the Lagrange interpolation on $\Xi$. Obviously, $P$ is
finite-rank since its range is a $\#\Xi$-dimensional subspace
of $\mathbb{F}[\bm{x}]$. Lagrange projectors are the standard
examples of ideal projectors.

It's well-known that an ideal projector can be characterized completely by the range of its dual projector, cf. \cite{MFS1916, Moller1977, deboorRon1991, MMM1993}.
Given a finite-rank linear projector $P$ on $\mathbb{F}[\bm{x}]$,
the kernel of $P$ is an ideal if and only if
the range of its dual projector is of the form
$$\bigoplus_{\bm{\xi}\in\Xi} \delta_{\bm{\xi}}\mathcal {Q}_{\bm{\xi}}(D)$$
with some finite point set $\Xi\subset \mathbb{F}^d$,  $D$-invariant finite-dimensional polynomial subspace
$\mathcal{Q}_{\bm{\xi}}\subset \mathbb{F}[\bm{x}]$ for
each $\bm{\xi}\in \Xi$.
$\delta_{\bm{\xi}}$ denotes the evaluation functional at the point $\bm{\xi}$,
and $\mathcal {Q}_{\bm{\xi}}(D)$ will be explained in the next section.

In the univariate case, for an integer $n$, there is only one $D$-invariant polynomial subspace of degree less than $n$, which implies
that
every univariate ideal projector can be viewed as a limiting case of Lagrange projector, cf. \cite{deboor2006}.
This prompted C. de Boor to define \emph{Hermite projector} as the pointwise limit of Lagrange projectors and
pose the following conjecture in \cite{dBo2005}. Indeed, this conjecture had been raised in \cite{BR1990}
with certain restriction.

\hspace{1pt}\\
\textbf{C. de Boor's conjecture} \emph{A finite-rank linear projector on
$\mathbb{C}[\bm{x}]$ is an ideal  projector
if and only if it is the pointwise
limit of Lagrange projectors}.
\\

B. Shekhtman \cite{BS2006} constructed a counterexample to this conjecture for every $d\geq 3$.
In the same paper, B. Shekthman also
showed that the conjecture is true for bivariate complex projectors
with the help of  Fogarty Theorem (see \cite{Foga1968}). Later,
using the fact that any pair of commuting matrices can be approximated by pairs of diagonalizable commuting
complex matrices(see \cite{TT1955, GR1992}), C. de Boor and B. Shekhtman \cite{deBoorShe2008} reproved the
same result. In addition, by Theorem 8 of \cite{GRSB2000},  \cite{deBoorShe2008} also proved that certain low-rank multivariate ideal projectors are
the limit of Lagrange projectors.
Specifically, B. Shekhtman \cite{BS2008} completely analyzed the
bivariate ideal projectors which are onto the space of polynomials
of degree less than $n$ over real or complex field, and verified the
conjecture in this particular case.

Since for every $d\geq 3$, there exist ideal projectors that are not
the pointwise limits of Lagrange projectors, the question that what
kind of ideal projectors can be perturbed as Lagrange projectors
lies ahead. For this purpose, B. Shekhtman \cite{BS2009} theoretically presented
a symbolic algorithm which can determine whether an ideal projector
is the limit of Lagrange projectors or not. However, as mentioned by this paper,
such a
method isn't yet feasible in practice.

In the converse case,
one wonders how to generate a sequence of Lagrange projectors practically
such that this sequence of Lagrange projectors converges pointwise to a given Hermite projector. More deeply,
B. Shekthman  raised the following question in \cite{She2009}.

\hspace{1pt}\\
\textbf{Problem}  \emph{
Let $P$ be an Hermite projector, and $P_h$ a sequence of Lagrange projectors such that $P_h$ converges pointwise to $P$ as $h$ tends to zero.
Then, what's
the relationship between the trajectories of the points in varieties of $\mathrm{ker} P_h$ and $P$}?
\\

In this paper, we deal with two classes of $D$-invariant subspaces. The first one is classical, which
has been investigated by many literatures such as \cite{BR1990, SX1995,  Lor00, GS00:3, HH2000}.
The second  is some special type, which is inspired by examples in \cite{BR1990, dBo2005, She2009}.
We construct a group of interpolation point sets corresponding to each class, and
establish the relationship between evaluation functionals induced by the point sets and derivative functionals
related to the corresponding $D$-invariant subspaces.
Based on these results,  we generate a sequence of Lagrange projectors which converges pointwise to the ideal projector
with related $D$-invariant subspaces belonging to either the first or the second class.
Equivalently,
we constructively verify  C. de Boor's conjecture for the
ideal projector of this type.

The remainder of this paper is organized as follows. The next section
is devoted as a preparation for the paper. Section 3 and 4 discuss the first class of $D$-invariant subspaces and the second, respectively.
Finally,  main theorem of this paper is given in Section 5.

\section{Preliminaries}\label{pre}

In this section, we will firstly introduce some notations and review some
basic facts about ideal projector. For more details, we refer
the reader to \cite{dBo2005, She2009}.

Throughout the paper, we use $\mathbb{N}$ to stand for the monoid of
nonnegative integers and boldface type for tuples with their entries
denoted by the same letter with subscripts, for example,
$\bm{\alpha}=(\alpha_1,\ldots, \alpha_d)$.
For arbitrary $\bm{\alpha}\in \mathbb{N}^d$,  we define $\|\bm{\alpha}\|_1=\alpha_1+\ldots+\alpha_d$ and $\bm{\alpha}!=\alpha_1!\ldots \alpha_d !$.

For arbitrary tuples
$\bm{\alpha}$, $\bm{\beta}\in\mathbb{N}^d$, we write $\bm{\beta}\leq \bm{\alpha}$ if
$\bm{\alpha}-\bm{\beta}$ has only nonnegative entries. In other words, $\leq$
is the usual product order on $\mathbb{N}^d$.
A subset $\mathfrak{\Delta}\subset \mathbb{N}^d$ is called a
\emph{lower set} (alternatively, \emph{down set}, \emph{order
ideal}, etc.) if it is \emph{closed} under $\leq$, that is, $\bm{\alpha}
\in\mathfrak{\Delta}$ implies $\bm{\beta} \in\mathfrak{\Delta}$ for all
$\bm{\beta}\leq \bm{\alpha}$.

Let $P$ be a finite-rank ideal projector on $\mathbb{F}[\bm{x}]$.
The range and kernel of $P$ will be described as
\begin{align*}
        \mathrm{ran}P:=&P(\mathbb{F}[\bm{x}])=\{g \in \mathbb{F}[\bm{x}]: g=Pf \mbox{~for~some~} f\in \mathbb{F}[\bm{x}]\},\\
        \mathrm{ker}P:=&\mathrm{ran}(I-P)=\{g \in \mathbb{F}[\bm{x}]: Pg=0\},
\end{align*}
where $\mathrm{ran}P$ forms a finite-dimensional subspace and
$\mathrm{ker}P$ a zero-dimensional ideal of $\mathbb{F}[\bm{x}]$. It
is  obvious that the ideal $\mathrm{ker}P$ complements
the subspace $\mathrm{ran} P$ , i.e., $\mathrm{ker}P\oplus
\mathrm{ran}P=\mathbb{F}[\bm{x}]$.

Furthermore, as an infinite dimensional $\mathbb{F}$-vector space,
$\mathbb{F}[\bm{x}]$ has an algebraic dual $\mathbb{F}[\bm{x}]'$. An
ideal projector $P$ on $\mathbb{F}[\bm{x}]$ also has a dual
projector $P'$ on $\mathbb{F}[\bm{x}]'$, and
$$\mathrm{ran}P'={(\mathrm{ker} P)}^{\perp}:=\{\lambda \in \mathbb{F}[\bm{x}]': \mathrm{ker}P\subset \mathrm{ker}\lambda \}.$$
In fact, $\mathrm{ran}P'$ is the set of interpolation conditions
matched by $P$. It is easy to see that the maximum number of
linearly independent functionals in $\mathrm{ran}P'$ equals the dimension of
$\mathrm{ran}P$, namely $\dim\mathrm{ran} P=\dim\mathrm{ran}P'$. In addition,
if $\bm{q}=(q_1,\ldots,q_s)$
is an $\mathbb{F}$-basis for $\mathrm{ran}P$ and
$\bm{\lambda}=(\lambda_1,\ldots,\lambda_s)$ an $\mathbb{F}$-basis
for $\mathrm{ran}P'$, then their  Gram matrix
$$\bm{\lambda}^T \bm{q}:=(\lambda_i q_j)_{1\leq i, j\leq s}$$ is invertible.

For $\bm{\alpha}\in \mathbb{N}^d $, we write $\bm{x}^{\bm{\alpha}}$
for the monomial $x_1^{\alpha_1}\ldots x_d^{\alpha_d}$, especially define $\bm{x}^{\bm{\alpha}}=1$ for $\bm{\alpha}=(0,\ldots , 0)$.
 Thus, a
polynomial $p$ in $\mathbb{F}[\bm{x}]$ can be expressed as
\begin{equation}\label{pform}
p=\sum\limits_{\bm{\alpha}\in \mathbb{N}^d}\widehat{p}({\bm{\alpha}})
{\bm{x}}^{{\bm{\alpha}}}
\end{equation}
with $\widehat{p}({\bm{\alpha}})\in \mathbb{F}$ nonzero.
For a polynomial
$p(\bm{x})$ given as in (\ref{pform}), we write
\begin{align*}
p(D): \mathbb{F}[\bm{x}]&\rightarrow \mathbb{F}[\bm{x}]\\
 f&\mapsto \sum\limits_{\bm{\alpha}\in \mathbb{N}^d}\widehat{p}({\bm{\alpha}})
\frac{\partial
^{\|\bm{\alpha}\|_1}f}{\partial x_1^{\alpha_1}\cdots\partial
x_d^{\alpha_d}}
\end{align*}
for the associated differential operator. For a finite-dimensional polynomial subspace
$\mathcal {Q}\subset \mathbb{F}[\bm{x}]$,
we define
$$\mathcal {Q}(D):=\mathrm{span}_{\mathbb{F}}\{q(D): q\in \mathcal {Q}\}.$$
We call a polynomial subspace $\mathcal {Q}$
a
\emph{$D$-invariant polynomial subspace}
if it's closed under differentiation, i.e., for every
$ q\in \mathcal {Q}$,
$\frac{\partial q}{\partial x_i}\in \mathcal {Q}$ for all $1\leq i\leq d$.

Finally, let us recall some facts about combinatorics.
For $r, m \in \mathbb{N}$ with $0\leq r\leq m$,  $m \choose r$ is the  binomial coefficient, i.e., ${m \choose r}=\frac{m!}{r!(m-r)!}$.
For $\bm{\alpha}=(\alpha_1, \ldots, \alpha_d)\in \mathbb{N}^d$, ${\|\bm{\alpha}\|_1 \choose \alpha_1, \ldots, \alpha_d}$ signifies multinomial coefficient, that is,
$${\|\bm{\alpha}\|_1 \choose \alpha_1, \alpha_2, \ldots, \alpha_d}=\frac{\|\bm{\alpha}\|_1!}{\bm{\alpha}!}.$$

\begin{lem}\textup{\cite[p. 90]{Combina2001}} \label{L1}
Let $n, m$ be arbitrary nonnegative integers satisfying  $n\leq m$. Then
$$
\sum\limits_{r=0}^{m}(-1)^{m-r} {m \choose r}r^n=\left\{
                                                \begin{array}{ll}
                                                  m!, &n=m; \\
                                                  0, & 0 \leq n<m.
                                                \end{array}
                                              \right.
$$
\end{lem}

\section{The first class of $D$-invariant subspaces}

In this section, we are concerned with the classical $D$-invariant subspaces, which only involve directional derivatives.
Based on the study of \cite{BR1990, dBo2005, GS00:3}, we describe that how the interpolation points are arranged in the straight directions,
this sequence of evaluation functionals can
converge to the corresponding directional derivative functional.

Let $\bm{\alpha}=(\alpha_1, \ldots, \alpha_d)\in \mathbb{N}^d$ and $\bm{\rho}=(\bm{\rho}_{1}, \bm{\rho}_{2}, \ldots, \bm{\rho}_{d})\in (\mathbb{F}^d)^d$ in which
$$\bm{\rho}_{1}=({\rho}_{1, 1},
\ldots, {\rho}_{1, d}), \bm{\rho}_{2}=({\rho}_{2, 1},
\ldots, {\rho}_{2, d}), \ldots, \bm{\rho}_{d}=({\rho}_{d, 1},
\ldots, {\rho}_{d, d})\in \mathbb{F}^d$$
are $\mathbb{F}$-linearly independent unit vectors. Then we define a differential operator
$D_{\bm{\rho}}^{\bm{\alpha}}: \mathbb{F}[\bm{x}]\rightarrow \mathbb{F}[\bm{x}]$ by the formula
\begin{align*}
D_{\bm{\rho}}^{\bm{\alpha}}f=
\frac{\partial
^{\|\bm{\alpha}\|_1}f}{\partial \bm{\rho}_1^{\alpha_1}\cdots\partial
\bm{\rho}_d^{\alpha_d}}.
\end{align*}
Furthermore,
the next lemma tells that for what type of polynomial, its associated differential operator is exactly $D_{\bm{\rho}}^{\bm{\alpha}}$.
\begin{lem}\label{D1lem}
Let $\bm{\alpha}=(\alpha_1, \ldots, \alpha_d)\in \mathbb{N}^d$, $\bm{\rho}=(\bm{\rho}_{1}, \bm{\rho}_{2}, \ldots, \bm{\rho}_{d})\in (\mathbb{F}^d)^d$ be as above,
and let
$$
p=\prod_{i=1}^d(\bm{\rho}_{i}\cdot \bm{x})^{\alpha_i},
$$
where $\cdot$ denotes the Euclidean product.
Then
$p(D)=D_{\bm{\rho}}^{\bm{\alpha}}.$
\end{lem}
\begin{pf}
The proof can be easily completed by direct computation.\qed
\end{pf}

\begin{prop}\label{D1prop1}
Given a lower set $\mathfrak{\Delta} \subset\mathbb{N}^d$ and $\bm{\rho}=(\bm{\rho}_{1}, \bm{\rho}_{2}, \ldots, \bm{\rho}_{d})\in (\mathbb{F}^d)^d$ as above,
then the subspace
$$
\mathcal {Q}=\mathrm{span}_{\mathbb{F}}\{\prod_{i=1}^d(\bm{\rho}_{i}\cdot \bm{x})^{\alpha_i} :\bm{\alpha}\in
\mathfrak{\Delta}\}\subset \mathbb{F}[\bm{x}].
$$
is $D$-invariant.
\end{prop}
\begin{pf}
This proposition is the immediate consequence of the differential calculus. \qed
\end{pf}

The following proposition makes a
connection between difference quotient and directional derivative of
multivariate polynomial, which is a variant of  Proposition 9.2 of \cite{BR1990} and  Proposition 7.2 of \cite{dBo2005}. We
also provide a simple proof, for completeness.

\begin{prop}\label{D1prop2}
Let $\bm{\rho}=(\bm{\rho}_{1}, \bm{\rho}_{2}, \ldots, \bm{\rho}_{d})\in (\mathbb{F}^d)^d$ be as above,
$h$ a non-zero number in $\mathbb{F}$, and let
$\bm{\alpha}\in \mathbb{N}^d$, $\bm{\xi}\in \mathbb{F}^d$.
Then for arbitrary $p\in \mathbb{F}[\bm{x}]$,
\begin{equation}\label{mianformulas1}
\begin{aligned}
\frac{1}{{h}^{\|\bm{\alpha}\|_1}}\sum_{\bm{0}\leq \bm{\beta}\leq \bm{\alpha}}(-1)^{\|\bm{\alpha}-\bm{\beta}\|_1}
{{\bm{\alpha}} \choose {\bm{\beta}}}p({\bm{\xi}+ h\sum_{i=1}^d  \beta_i\bm{\rho}_i})
=(D_{\bm{\rho}}^{\bm{\alpha}}p)(\bm{\xi})+O(h),
\end{aligned}
\end{equation}
where ${{\bm{\alpha}} \choose {\bm{\beta}}}=\prod_{i=1}^d {\alpha_i \choose \beta_i}$,  and the remainder $O(h)$ is a polynomial in $h$.
\end{prop}
\begin{pf}
For convenience, let
 $$\varphi=\sum_{\bm{0}\leq \bm{\beta}\leq \bm{\alpha}}(-1)^{\|\bm{\alpha}-\bm{\beta}\|_1}
{{\bm{\alpha}} \choose {\bm{\beta}}} p(\bm{\xi}+ h\sum_{k=1}^d  \beta_k\bm{\rho}_k).$$
Then we can obtain that for $m=0,\ldots, \|\bm{\alpha}\|_1$,
\begin{align*}
\frac{d^m \varphi(0)}{d h^m}=&
\sum_{\|\bm{\gamma}_1 \|_1+\ldots+\|\bm{\gamma}_d \|_1=m }
\rho_{1, 1}^{\gamma_{1,1}}\ldots \rho_{d, 1}^{\gamma_{1,d}}\ldots
\rho_{1, d}^{\gamma_{d,1}}\ldots \rho_{d,
d}^{\gamma_{d,d}}\frac{\partial^m p(\bm{\xi})}{\partial x_1^{\|\bm{\gamma}_1 \|_1
}\ldots
\partial x_d^{\|\bm{\gamma}_d \|_1}
 }\\
&\left(
                                                               \begin{array}{c}
                                                                 m\\
                                                               \gamma_{1, 1}\ldots   \gamma_{1, d},\ldots,  \gamma_{d, 1},  \ldots,  \gamma_{d, d}\\
                                                               \end{array}
                                                             \right)
\prod_{j=1}^d(\sum_{\beta_j=0}^{\alpha_j}
(-1)^{\alpha_j-\beta_j}{\alpha_j \choose \beta_j}\beta_j^{\sum_{i=1}^d \gamma_{i, j}})
\end{align*}
with $\bm{\gamma}_i=(\gamma_{i, 1}, \ldots, \gamma_{i, d})\in \mathbb{N}^d$, $1\leq i\leq d$.

In the following, we will simplify the right-hand side of the above equality.
There are two cases which must be examined.

Case 1:  $0\leq
m\leq \|\bm{\alpha}\|_1-1$.

In this case, there must exist some $1\leq j\leq d$ such
that
$\sum_{i=1}^d\gamma_{i, j}<\alpha_j$. By Lemma 1, it follows for such $j$,
$$\sum_{\beta_j=0}^{\alpha_j}
(-1)^{\alpha_j-\beta_j}{\alpha_j \choose \beta_j}\beta_j^{\sum_{i=1}^d\gamma_{i, j}}=0.$$
Thus, for all $0\leq m\leq \|\bm{\alpha}\|_1-1$,
$\frac{d^m \varphi(0)}{d
 h^m}=0$.

Case 2: $m=\|\bm{\alpha}\|_1$.

In this case, if there exists some $1\leq j\leq d$ such that $\sum_{i=1}^d\gamma_{i, j}>\alpha_j$, then
there must exist another $1\leq j'\leq d$,
such that $\sum_{i=1}^d\gamma_{i, j'}<\alpha_{j'}$.
Hence, when $m=\|\bm{\alpha}\|_1$, $$\prod_{j=1}^d(\sum_{\beta_j=0}^{\alpha_j}
(-1)^{\alpha_j-\beta_j}{\alpha_j \choose \beta_j}\beta_j^{\sum_{i=1}^d\gamma_{i, j}})\neq 0$$ if and only if for all $1\leq j\leq d$,
$\sum_{i=1}^d\gamma_{i, j}=\alpha_j.$
Combining this with the fact
 $$\sum_{\beta_j=0}^{\alpha_j}
(-1)^{\alpha_j-\beta_j}{\alpha_j \choose \beta_j}\beta_j^{\alpha_j}=\alpha_j!, \quad 1\leq j\leq d,$$
we can conclude that
\begin{align*}
\frac{d^{\|\bm{\alpha}\|_1}\varphi(0)}{d
h^{\|\bm{\alpha}\|_1}}=
\|\bm{\alpha}\|_1!&\sum_{\gamma_{1,1}+\ldots+\gamma_{d,
1}=\alpha_1}\ldots  \sum_{\gamma_{1,d}+\ldots+ \gamma_{d, d}=\alpha_d}\rho_{1, 1}^{\gamma_{1,1}}\ldots \rho_{d, 1}^{\gamma_{1,d}}\ldots
\rho_{1, d}^{\gamma_{d,1}}\ldots \rho_{d,
d}^{\gamma_{d,d}}\\
&
\left(
                                                               \begin{array}{c}
                                                                 \alpha_1\\
                                                                 \gamma_{1,1},\ldots, \gamma_{d, 1}\\
                                                               \end{array}
                                                             \right)\ldots\left(
                                                               \begin{array}{c}
                                                                 \alpha_d\\
                                                                 \gamma_{1, d},\ldots, \gamma_{d,d}\\
                                                               \end{array}
                                                             \right)\frac{\partial^{\|\bm{\alpha}\|_1} p(\bm{\xi})}{\partial x_1^{\|\bm{\gamma}_1\|_1
}\ldots
\partial x_d^{\|\bm{\gamma_d\|_1
}}}.
\end{align*}
That is,
$$\frac{d^{\|\bm{\alpha}\|_1}\varphi(0)}{d
h^{\|\bm{\alpha}\|_1}}
={\|\bm{\alpha}\|_1} !(D_{\bm{\rho}}^{\bm{\alpha}}p)(\bm{\xi}).$$

In sum, we have deduced that
$$\frac{d^m  \varphi(0)}{d h^m}=\left\{
    \begin{array}{ll}
      0, &   0\leq m \leq \|\bm{\alpha}\|_1-1; \\
    \|\bm{\alpha}\|_1!(D_{\bm{\rho}}^{\bm{\alpha}}p)(\bm{\xi}), &m=\|\bm{\alpha}\|_1,
    \end{array}
  \right.
$$
which implies that
$$\varphi=h^{\|\bm{\alpha}\|_1}(D_{\bm{\rho}}^{\bm{\alpha}}p)(\bm{\xi})+O(h^{\|\bm{\alpha}\|_1+1}).$$
Equality (\ref{mianformulas1}) follows directly from the above equality.
\qed
\end{pf}

\section{The second class of $D$-invariant subspaces}

In this section, we  present a new class of $D$-invariant subspaces, which are spanned by  polynomials with special structure.
We also provide a group of interpolation point
sets corresponding to this class of $D$-invariant subspaces,
and
discuss that how some interpolation points
coalesce in the non-straight directions, such sequence of  evaluation functionals can
converge to the differential functional about this class of $D$-invariant subspaces.

Before introducing the new class of $D$-invariant subspaces,
we need to settle some notations used throughout this section.

Let
$\bm{a}=(a_0, a_1, \ldots, a_n)$ with $n\geq 1$ be an $n+1$-tuple of positive integers satisfying
\begin{equation}\label{2ajcondition}
a_0=1 \mbox{~and~} a_1>\ldots>a_n\geq 2,
\end{equation}
and let
$$
\bm{c}_1=(c_{1, 0}, c_{1, 1}, \ldots, c_{1, n}),\ldots, \bm{c}_d=(c_{d, 0}, c_{d, 1}, \ldots, c_{d, n})\in \mathbb{F}^{n+1}
$$
satisfying that
$c_{1, 0}, c_{2, 0}, \ldots, c_{d, 0}$ aren't all zero.
Then it's clear that
\begin{equation}\label{2map}
\tau(\bm{\gamma}_{1}, \ldots, \bm{\gamma}_{d})=
 \sum_{j=0}^n a_j\sum_{i=1}^d \gamma_{i, j}
 \end{equation}
with $\bm{\gamma}_{i}=(\gamma_{i, 0}, \gamma_{i, 1}, \ldots, \gamma_{i, n})\in \mathbb{N}^{n+1}$ defines a map
$\tau:(\mathbb{N}^{n+1})^{d} \rightarrow \mathbb{N}$.

\begin{prop}\label{2D2prop1}
Let $\bm{a}=(a_0, a_1, \ldots, a_n)$, $\bm{c}_{i}=(c_{i, 0}, c_{i, 1}, \ldots, c_{i, n})$, $1\leq i\leq d$ 
and the map $\tau$ be as
above. Let $q_{n, m}$, $m=0, 1, \ldots, a_1$ be polynomials defined by
$$
q_{n, m}=\sum_{
\tau(\bm{\gamma}_{1}, \ldots, \bm{\gamma}_{d})=m}
\frac{\bm{c}_1^{\bm{\gamma}_1}\ldots \bm{c}_d^{\bm{\gamma}_d}}{\bm{\gamma}_1 !\ldots \bm{\gamma}_d !}
x_1^{\| \bm{\gamma}_1\|_1 }\ldots
x_d^{\| \bm{\gamma}_d\|_1 },
$$
with $\bm{\gamma}_{i}=(\gamma_{i, 0}, \gamma_{i, 1}, \ldots, \gamma_{i, n})\in \mathbb{N}^{n+1}$,
$\bm{c}_i^{\bm{\gamma}_i}=\prod_{j=0}^n c_{i, j}^{\gamma_{i, j}}$,  and let $\mathcal {Q}$ be a polynomial subspace defined  by
$$
\mathcal {Q}=\mathrm{span}_{\mathbb{F}}\{q_{n, m}: 0\leq m\leq a_1\}.
$$
Then the following hold:
\begin{enumerate}
  \item[\textup{(i)}] $\mathrm{dim} \mathcal {Q}=a_1+1$,
  \item[\textup{(ii)}] $\mathcal {Q}$ is a $D$-invariant polynomial subspace.
\end{enumerate}
\end{prop}
\begin{pf}
To prove \textup{(i)},
we assume that there exist $k_0, k_1, \ldots, k_{a_1-1}, k_{a_1}\in \mathbb{F}$ such that
$$k_0 q_{n, 0}+k_{1} q_{n, 1}+\ldots+ k_{a_1-1} q_{n, a_1-1}+k_{a_1} q_{n, a_1}=0.$$
Since there exists some $1\leq i_0\leq d$ such that $c_{i_0, 0}\neq 0$, then $x_{i_0}^{m}$ must belong to the support of $q_{n, m}$.
Using this together with the definition of $\tau$, we conclude that the degree of $q_{n, m}$ is $m$.

Specifically, $x_{i_0}^{a_1}$  belongs to the support of $q_{n, a_1}$,  while for all
$0\leq m\leq a_1-1$,
$x_{i_0}^{a_1}$ can't belong to the support of $q_{n, m}$. Therefore $k_{a_1}=0$, which implies that
$$k_0 q_{n, 0}+k_{1} q_{n, 1}+\ldots+ k_{a_1-1} q_{n, a_1-1}=0.$$

Arguing for $m=a_1-1, a_1-2, \ldots, 1, 0$ as for $m=a_1$,  we get
$k_{a_1-1}=0, k_{a_1-2}=0,  \ldots, k_{1}=0, k_0=0$, successively.
That is to say, the family of polynomials $q_{n, m}, 0\leq m\leq a_1$ is $\mathbb{F}$-linearly independent.
So $\mathrm{dim} \mathcal {Q}=a_1+1$.

To prove \textup{(ii)}, we firstly notice that if $1\leq m< a_n$,  then
$\tau(\bm{\gamma}_{1}, \ldots, \bm{\gamma}_{d})=m$ implies that $\gamma_{i, j}=0$ with $1\leq i\leq d, 1\leq j\leq n$. Thus,  for all $1\leq m< a_n$,
\begin{align*}
q_{n, m}&=\sum_{\gamma_{1, 0}+\ldots+ \gamma_{d, 0}=m}\frac{c_{1, 0}^{\gamma_{1, 0}}\ldots c_{d, 0}^{\gamma_{d, 0}}}{\gamma_{1, 0}!\ldots \gamma_{d, 0}!}x_1^{\gamma_{1, 0}}\ldots x_d^{\gamma_{d, 0}}\\
&=\frac{1}{m!}(c_{1, 0}x_1+c_{2, 0}x_2+\ldots+c_{d, 0}x_d)^m.
\end{align*}
Consequently,
\begin{equation}\label{2Dsimp}
\frac{\partial q_{n, m}}{\partial x_i}=c_{i, 0} q_{n, m-1}, \quad 1\leq m< a_n.
\end{equation}

Next, let $m$ be an arbitrary integer satisfying $a_n\leq m \leq a_1$, and $s$
the minimum integer  between $1$ and $n$ such that $m-a_s\geq 0$.
We claim
\begin{equation}\label{2conc3}
\frac{\partial q_{n, m}}{\partial x_i}=c_{i, 0}q_{n, m-1}+\sum_{j=s}^{n}c_{i, j}q_{n, m-a_j}, \quad a_n\leq m \leq a_1.
\end{equation}
This claim together with equality (\ref{2Dsimp})
immediately means that $\mathcal {Q}$ is $D$-invariant.

To prove our claim, we will use induction on the number $n$.
When $n=1$, our claim can be easily verified.
Now, assume that our claim is true for $n-1$. To prove that it holds for $n$,
let $k$ be the maximum nonnegative integer such that $k a_n \leq m$,
that is,
\begin{equation}\label{2kcondition}
k a_n \leq m, (k+1)a_n>m.
\end{equation}
and $k'$ the maximum nonnegative integer such that $k' a_n \leq m-1$.
From the definition of  $q_{n, m}$, we obtain that
\begin{align*}
q_{n, m}&=\sum_{l=0}^{k}
\big(\sum_
{\tau((\gamma_{1, 0}, \ldots ,\gamma_{1, n-1}), \ldots, (\gamma_{d, 0}, \ldots ,\gamma_{d, n-1}))=m-l a_n }\frac{c_{1, 0}^{\gamma_{1, 0}}\ldots c_{1, n-1}^{\gamma_{1, n-1}}\ldots c_{d, 0}^{\gamma_{d, 0}}\ldots c_{d, n-1}^{\gamma_{d, n-1}}}{
\gamma_{1, 0}!\ldots \gamma_{1, n-1}!\ldots \gamma_{d, 0}!\ldots \gamma_{d, n-1}!
}\\
&~~~~~x_1^{\gamma_{1, 0}+\ldots+\gamma_{1, n-1}}\ldots x_d^{\gamma_{d, 0}+\ldots+\gamma_{d, n-1}}\big)
\big(\sum_{\gamma_{1, n}+\ldots+ \gamma_{d, n}=l}
\frac{c_{1, n}^{\gamma_{1, n}}\ldots c_{d, n}^{\gamma_{d, n}}}{c_{1, n}! \ldots c_{d, n}!}x_1^{\gamma_{1, n}}\ldots x_d^{\gamma_{d, n}}\big).
\end{align*}
That is,
\begin{equation}\label{2tidui}
q_{n, m}=\sum_{l=0}^{k}\frac{1}{l !}(c_{1, n} x_1+\ldots+c_{d, n} x_d)^l q_{n-1, m-l a_n},
\end{equation}
which plays an important role in what follows.
At this point, we have two cases to consider.

Case 1: $a_n\leq m<a_{n-1}$.

In this case, $s=n$ and $0\leq m-l a_n<a_{n-1}$ for all $0\leq l\leq k$.
Hence, $$q_{n-1, m-l a_n}=\frac{1}{(m-l a_n)!}(c_{1, 0} x_1+\ldots+c_{d, 0} x_d)^{m-l a_n }, \quad 0\leq l\leq k.$$
By (\ref{2tidui}), we get
$$q_{n, m}=\sum_{l=0}^{k}\frac{1}{l !(m-l a_n)!}(c_{1, n} x_1+\ldots+c_{d, n} x_d)^l (c_{1, 0} x_1+\ldots+c_{d, 0} x_d)^{m-l a_n }.$$
So we obtain that
\begin{align*}
&\frac{\partial q_{n, m}}{\partial x_i}\\
=&c_{i, 0}\sum_{l=0}^{k'}\frac{1}{l !(m-1-l a_n)!}(c_{1, n} x_1+\ldots+c_{d, n} x_d)^l (c_{1, 0} x_1+\ldots+c_{d, 0} x_d)^{m-1-l a_n }+\\
&c_{i, n}\sum_{l=1}^{k}\frac{1}{(l-1) !(m-l a_n)!}(c_{1, n} x_1+\ldots+c_{d, n} x_d)^{l-1} (c_{1, 0} x_1+\ldots+c_{d, 0} x_d)^{m-l a_n }\\
=&c_{i, 0}q_{n, m-1}+c_{i, n}q_{n, m-a_n}.
\end{align*}

Case 2: $a_{n-1}\leq m \leq a_1$.

In this case,  $0\leq s\leq n-1$.

By (\ref{2kcondition}) and (\ref{2ajcondition}),
it follows that $m-k a_n-a_{n-1}<0$. Thus,
there exists some $0\leq l_0\leq k-1$ such that
\begin{equation}\label{2l0condition}
m-l_0 a_n-a_{n-1}\geq 0, m-(l_0+1) a_n-a_{n-1}< 0.
\end{equation}
As a result, we have that $0\leq m-l a_n<a_{n-1}$ for all $l_0+1 \leq l\leq k$, which implies that
$$
q_{n-1, m-l a_n}=\frac{1}{ (m-l a_n)!}(c_{1, 0}x_1+\ldots+ c_{d, 0}x_d )^{m-l a_n}, \quad l_0+1 \leq l\leq k.
$$
Therefore,
\begin{equation}\label{2qiudao1}
\frac{\partial q_{n-1, m-l a_n}}{\partial x_i}=\left\{
  \begin{array}{ll}
    c_{i, 0} q_{n-1, m-1-l a_n},  & l_0+1 \leq l\leq k' \hbox{;} \\
    0, &   k'+1 \leq l\leq k \hbox{.}
  \end{array}
\right.
\end{equation}

For all $0\leq l\leq l_0$, we have $a_{n-1}\leq  m-l a_n\leq a_1$. Then our inductive hypothesis implies that
\begin{equation}\label{2qiudao2}
\frac{\partial q_{n-1, m-l a_n}}{\partial x_i}=c_{i, 0}q_{n-1, m-1-l a_n}+\sum_{j=s_l}^{n-1}c_{i, j}q_{n-1, m-l a_n-a_j}, \quad 0\leq l\leq l_0,
\end{equation}
where $s_{l}$ is the minimal integer  between $1$ and $n-1$ such that $m-l a_n-a_{s_l}\geq 0$.
It should be noticed that $s_{0}=s$.

By (\ref{2qiudao1}) and (\ref{2qiudao2}), we can deduce that
\begin{align*}
&\frac{\partial q_{n, m}}{\partial x_i}\\
=&c_{i, 0}(\sum_{l=0}^{k'}\frac{1}{l !}(c_{1, n} x_1+\ldots+c_{d, n} x_d)^l q_{n-1, m-1-la_n})+\\
&c_{i, n}
(\sum_{l=1}^{k}\frac{1}{(l-1) !}(c_{1, n} x_1+\ldots+c_{d, n} x_d)^{l-1}q_{n-1, m-l a_n})+\\
&\sum_{l=0}^{l_0}\frac{1}{l !}(c_{1, n} x_2+\ldots+c_{d, n} x_d)^l\sum_{j=s_l}^{n-1}c_{i, j}q_{n-1, m-l a_n-a_j}\\
=&c_{i, 0} q_{n, m-1}+c_{i, n}q_{n, m-a_n}+\sum_{l=0}^{l_0}\frac{1}{l !}(c_{1, n} x_2+\ldots+c_{d, n} x_d)^l\sum_{j=s_l}^{n-1}c_{i, j}q_{n-1, m-l a_n-a_j}.
\end{align*}

It remains to show the last row of the above equality and  the right-side hand of (\ref{2conc3}) are equal. More precisely,
for each $s_0\leq j\leq n-1$, let $k_j$ denote the maximum nonnegative integer  satisfying $k_j a_n\leq m-a_j$,
that is,
\begin{equation}\label{2kjcondition}
k_j a_n\leq m-a_j,  (k_j+1)a_n>m-a_j.
\end{equation}
From (\ref{2l0condition}), we know $k_{n-1}=l_0$.

Due to (\ref{2tidui}), we observe that
$$
q_{n, m-a_j}
=
\sum_{l=0}^{k_j}\frac{1}{l !}(c_{2, n} x_2+\ldots+c_{d, n} x_d)^{l}q_{n-1, m-a_j-l a_n},\quad s_0\leq j\leq n-1.
$$
Furthermore,
\begin{equation}\label{2you1}
q_{n, m-a_j}
=\sum_{l=0}^{l_0}\frac{1}{l !}(c_{1, n} x_2+\ldots+c_{d, n} x_d)^{l}q_{n-1, m-a_j-l a_n},\quad s_0\leq j\leq n-1.
\end{equation}
Here, we have used the fact that $q_{n-1, m-a_j-l a_n}=0$ for all $s_0\leq j\leq n-1$, $k_j+1\leq  l \leq l_0$.

Now, recall that  $s_{l}, 0\leq l\leq l_0$ is the minimal integer  between $1$ and $n-1$ such that $m-l a_n-a_{s_l}\geq 0$, then
\begin{equation}\label{2you3}
\sum_{j=s_0}^{n-1}c_{i, j}q_{n-1, m-a_j-l a_n}=\sum_{j=s_l}^{n-1}c_{i, j}q_{n-1, m-a_j-l a_n},  \quad 0\leq l\leq l_0.
\end{equation}

According to (\ref{2you1}), (\ref{2you3}) and the fact $s=s_0$, we find that
\begin{align*}
\sum_{j=s}^{n-1}c_{i, j}q_{n, m-a_j}
&=\sum_{j=s_0}^{n-1}c_{i, j} \sum_{l=0}^{l_0}\frac{1}{l !}(c_{1, n} x_2+\ldots+c_{d, n} x_d)^{l}q_{n-1, m-a_j-l a_n}\\
&=\sum_{l=0}^{l_0}\frac{1}{l !}(c_{1, n} x_2+\ldots+c_{d, n} x_d)^{l}(\sum_{j=s_l}^{n-1}c_{i, j}q_{n-1, m-a_j-l a_n}).
\end{align*}
Consequently,
$$\frac{\partial q_{n, m}}{\partial x_i}=c_{i, 0}q_{n, m-1}+\sum_{j=s}^{n}c_{i, j}q_{n, m-a_j}, \quad a_{n-1}\leq m\leq a_1.$$
\qed
\end{pf}

Next, we will give some simple examples of $D$-invariant subspaces as in Proposition \ref{2D2prop1}.

\begin{exmp}\label{e1}
Let $c_{i, j}=0, 1\leq i\leq d, 1\leq j\leq n$, then
$$\mathcal {Q}=\mathrm{span}_{\mathbb{F}}\{\frac{1}{m}(c_{1, 0}x_1+\ldots+c_{d, 0}x_d)^m: 0\leq m\leq a_1
\},
$$
which is spanned by  homogenous polynomials.  We call this case the trivial one. \qed
\end{exmp}
\begin{exmp}\label{e2}
Let $n=1$ and $c_{1, 0}=1$, $c_{i, 0}=0$ with $2\leq i\leq d$,
then $$\mathcal {Q}=\mathrm{span}_{\mathbb{F}}\{1, x_1, \frac{1}{2}x_1^2, \ldots, \frac{1}{a_1}x_1^{a_1}+c_{1, 1} x_1+c_{2, 1}x_2+\ldots+c_{d, 1}x_d\}.
$$

In this example, we should notice that when
$a_1=2$, $d=2$ and $c_{1, 1}=0, c_{2, 1}=1$, then $\mathcal {Q}=\mathrm{span}_{\mathbb{F}}\{1, x_1, \frac{1}{2}x_1^2+x_2\}$.
This is the example that had been discussed by some well-known papers,
for instance, \cite[p. 301]{BR1990}, \cite[p. 81]{dBo2005}, and Illustration 6.1.9 of \cite{She2009}.
\qed
\end{exmp}

\begin{exmp}\label{e3}
Let $n=1, \bm{a}=(1, 2), d=3, \bm{c}_1=(1, 0), \bm{c}_2=(1, 1), \bm{c}_3=(0, 1)$, then
$$\mathcal {Q}=\mathrm{span}_{\mathbb{F}}\{1, x_1+x_2,
\frac{1}{2}x_1^2+x_1 x_2+\frac{1}{2}x_2^2+x_2+x_3\},$$ which will be considered in the next section.
\qed
\end{exmp}

The next proposition not only
gives us a type of interpolation point sets corresponding to the $D$-invariant subspaces as in Proposition \ref{2D2prop1},
but also reveals  their relationship.
\begin{prop}\label{2D2prop2}
Let $\bm{a}=(a_0, a_1, \ldots, a_n)$, $\bm{c}_{i}=(c_{i, 0}, c_{i, 1}, \ldots, c_{i, n})$ with $1\leq i\leq d$, and let
$q_{n, m}$ with $0\leq m\leq a_1$ be as above, $h$ a non-zero number in $\mathbb{F}$.  Then for arbitrary $p\in \mathbb{F}[\bm{x}]$ and $\bm{\xi}\in \mathbb{F}^d$,
\begin{equation}\label{2mianformulas2}
\begin{aligned}
&\frac{1}{{m!h^{m}}}\sum_{0\leq r\leq m}(-1)^{m-r} {m \choose r} p({\bm{\xi}+(\sum_{j=0}^n c_{1, j}(rh)^{a_j},  \ldots, \sum_{j=0}^n c_{d, j}(rh)^{a_j} )})
=\\&(q_{n, m}(D)p)(\bm{\xi})+O(h).
\end{aligned}
\end{equation}
\end{prop}
\begin{pf}
Applying Taylor Formulas, we obtain
\begin{align*}
&p(\bm{\xi}+(\sum_{j=0}^n c_{1, j}(rh)^{a_j}, \ldots, \sum_{j=0}^n c_{d, j}(rh)^{a_j} ))
=\\
&\sum_{k=0}^{\infty}\sum_{\| \bm{\gamma}_1\|_{1}+\ldots+\|\bm{\gamma}_d\|_{1}=k}
(r h)^{ \sum_{1\leq i\leq d, 0\leq j\leq n} a_j \gamma_{i, j}}
\frac{\bm{c}_1^{\bm{\gamma}_1}\ldots \bm{c}_d^{\bm{\gamma}_d}}{\bm{\gamma}_1 !\ldots \bm{\gamma}_d !}
\frac{\partial ^{k}p(\bm{\xi})}{\partial x_1^{\| \bm{\gamma}_1\|_{1}}
\ldots \partial
x_d^{\| \bm{\gamma}_d\|_{1}}},
\end{align*}
where $\bm{c}_i^{\bm{\gamma}_i}=\prod_{j=0}^n c_{i, j}^{\gamma_{i, j}}$,
 $\bm{\gamma}_i=({\gamma}_{i, 0}, \ldots, {\gamma}_{i, n})\in \mathbb{N}^{n+1}, 1\leq i\leq d$.

By means of the map $\tau$ defined as in (\ref{2map}), the above equality can be rewritten as
\begin{align*}
&p(\bm{\xi}+(\sum_{j=0}^n c_{1, j}(rh)^{a_j},  \ldots, \sum_{j=0}^n c_{d, j}(rh)^{a_j} ))=\\&
\sum_{l=0}^{m}\sum_{\tau(\bm{\gamma}_{1}, \ldots, \bm{\gamma}_{d})=l}
(r h)^{l}\frac{\bm{c}_1^{\bm{\gamma}_1}\ldots \bm{c}_d^{\bm{\gamma}_d}}{\bm{\gamma}_1 !\ldots \bm{\gamma}_d !}
\frac{\partial ^{\|\bm{\gamma}_1\|_{1}+\ldots+\|\bm{\gamma}_d\|_{1}}p(\bm{\xi})}{\partial
 x_1^{\| \bm{\gamma}_1\|_{1}}\ldots \partial
x_d^{\| \bm{\gamma}_d\|_{1}}}
+O(h^{m+1}),
\end{align*}
where $O(h^{m+1})$ is a polynomial in $h$.

Finally, by  Lemma \ref{L1},
we can conclude that
\begin{align*}
&\sum_{0\leq r\leq m}(-1)^{m-r}{m \choose r}
p(\bm{\xi}+(\sum_{j=0}^n c_{1, j}(rh)^{a_j},  \ldots, \sum_{j=0}^n c_{d, j}(rh)^{a_j} ))=
\\
&m! h^{m}
\sum_{\tau(\bm{\gamma}_{1}, \ldots, \bm{\gamma}_{d})=m}\frac{\bm{c}_1^{\bm{\gamma}_1}\ldots \bm{c}_d^{\bm{\gamma}_d}}{\bm{\gamma}_1 !\ldots \bm{\gamma}_d !}
\frac{\partial ^{\|\bm{\gamma}_1\|_{1}+\ldots+\|\bm{\gamma}_d\|_{1}}p(\bm{\xi})}{\partial
 x_1^{\| \bm{\gamma}_1\|_{1}}\ldots \partial
x_d^{\| \bm{\gamma}_d\|_{1}}}
+O(h^{m+1}),
\end{align*}
which leads to the proposition immediately.
\qed
\end{pf}

\section{Main theorem}

In this section, we consider a particular type of ideal projectors associated with the above two classes of $D$-invariant subspaces, and
constructively prove that C. de Boor's
conjecture is true for ideal projectors of this type.

First of all, we also introduce the notation that will be adopted in  main theorem.

Let  $\bm{\xi}^{(1)},\ldots,\bm{\xi}^{(\mu)}, \bm{\xi}^{(\mu+1)}, \ldots, \bm{\xi}^{(\mu+\nu)}\in \mathbb{F}^d$
be distinct points.
For each $1\leq k\leq \mu$, let $\mathfrak{\Delta}^{(k)} \subset\mathbb{N}^d$ be a lower set and
  $\bm{\rho}^{(k)}=(\bm{\rho}^{(k)}_1, \ldots, \bm{\rho}^{(k)}_d)\in (\mathbb{F}^d)^d$ be as in Section 3.
Likewise, for each $1\leq l\leq \nu$, let $\bm{a}^{(l)}=(a_0^{{(l)}}, a_1^{{(l)}}, \ldots, a_{n(l)}^{{(l)}})$,
  $\bm{c}_{i}^{(l)}=({c}_{i, 0}^{(l)}$, ${c}_{i, 1}^{(l)}, \ldots, {c}_{i, n(l)}^{(l)})$, $1\leq i\leq d$, and
$q_{n(l), m}^{(l)}$, $0\leq m\leq a_1^{(l)}$ be as in Section 4.

\begin{thm}\label{mianthm}
With the notation above,
let $P$ be an ideal projector with
\begin{align*}
\mathrm{ran} P'=\mathrm{span}_{\mathbb{F}}\big\{&\delta_{\bm{\xi}^{(k)}}D_{\bm{\rho}^{(k)}}^{\bm{\alpha}}, \delta_{\bm{\xi}^{(\mu+l)}}q_{n(l), m}^{(l)}(D)
:\\& \bm{\alpha}\in \mathfrak{\Delta}^{(k)}, 1\leq k\leq \mu, 0\leq m\leq a_1^{(l)},  1\leq l\leq \nu \big\},
\end{align*}
and let $P_h$ be a Lagrange projector with
\begin{align*}
\mathrm{ran} P'_h=
\mathrm{span}_{\mathbb{F}}\big\{&
\delta_{\bm{\xi}^{(k)}+ h\sum_{i=1}^d  \alpha_i\bm{\rho}_i^{(k)}},
\delta_{\bm{\xi}^{(\mu+l)}+\bm{\phi}^{(l)}(m h)}:\\
&\bm{\alpha}\in \mathfrak{\Delta}^{(k)}, 1\leq k\leq \mu, 0\leq m\leq a_1^{(l)}, 1\leq l\leq \nu \big\},
\end{align*}
where $h \in \mathbb{F}\setminus \{0\}$ and
$$\bm{\phi}^{(l)}(mh)=(\sum_{j=0}^{n(l)} c_{1, j}^{(l)}(m h)^{a_j^{(l)}}, \sum_{j=0}^{n(l)} c_{2, j}^{(l)}(m h)^{a_j^{(l)}}, \ldots, \sum_{j=0}^{n(l)} c_{d, j}^{(l)}{(m h)}^{a_j^{(l)}}).$$
Then the following statements hold:
\begin{enumerate}
  \item[\textup{(i)}] There exists a positive $\eta\in \mathbb{F}$ such that
   $$
   \mathrm{ran}P_h=\mathrm{ran}P, \quad \forall 0<|h|<\eta.
   $$
  \item[\textup{(ii)}] $P$ is the pointwise limit of $P_h, 0<|h|<\eta,$ as $h$ tends to zero.
\end{enumerate}
\end{thm}

\begin{pf}
Firstly, one can easily verify that
\begin{align*}
\bm{\lambda}=\big(&\delta_{\bm{\xi}^{(k)}}D_{\bm{\rho}^{(k)}}^{\bm{\alpha}}, \delta_{\bm{\xi}^{(\mu+l)}}q_{n(l), m}^{(l)}(D)
:\\& \bm{\alpha}\in \mathfrak{\Delta}^{(k)}, 1\leq k\leq \mu, 0\leq m\leq a_1^{(l)},  1\leq l\leq \nu\big)\in
{(\mathbb{F}[\bm{x}]')}^s
\end{align*}
and
\begin{align*}
\bm{\lambda}_h=\big(&\delta_{\bm{\xi}^{(k)}+ h\sum_{i=1}^d  \alpha_i\bm{\rho}_i^{(k)}},
\delta_{\bm{\xi}^{(\mu+l)}+\bm{\phi}^{(l)}(m h)}:\\
&\bm{\alpha}\in \mathfrak{\Delta}^{(k)}, 1\leq k\leq \mu, 0\leq m\leq a_1^{(l)}, 1\leq l\leq \nu \big)\in
{(\mathbb{F}[\bm{x}]')}^s
\end{align*}
form $\mathbb{F}$-bases for $\mathrm{ran} P'$ and $\mathrm{ran} P'_h$ respectively, where $$s=\sum_{k=1}^\mu \#\mathfrak{\Delta}^{(k)}+\sum_{k=1}^\nu a_1^{(k)}+\nu.$$
Provided that the entries of $\bm{\lambda}$ and
$\bm{\lambda}_h$ are arranged in the same order, namely for arbitrary fixed
 $1\leq k\leq \mu$, $\bm{\alpha}\in \mathfrak{\Delta}^{(k)}$,
the corresponding entries of $\bm{\lambda}$ and $\bm{\lambda}_h$ are in the same position,
the same as  for arbitrary fixed  $1\leq l\leq \nu$, $0\leq m\leq a_1^{(l)}$.  We
denote $\bm{\lambda}=(\lambda_1, \ldots, \lambda_s)$ and $\bm{\lambda}_h=(\lambda_{h, 1}, \ldots, \lambda_{h, s})$, respectively.

Let $\bm{q}=(q_1,q_2,\ldots,q_s)$ be an $\mathbb{F}$-basis for $\mathrm{ran}P$, For convenience, we introduce two $s \times s$ matrices
$$
\bm{\lambda}^T \bm{q}:=(\lambda_i q_j)_{1\leq i, j\leq s}, \quad \bm{\lambda}^T_h \bm{q}:=(\lambda_{h, i} q_j)_{1\leq i, j\leq s}.
$$
and for arbitrary $f \in \mathbb{F}[\bm{x}]$, vectors
$$
\bm{\lambda}^T f:=(\lambda_i f)_{1\leq i\leq s}, \quad \bm{\lambda}^T_h f:=(\lambda_{h, i} f)_{1\leq i\leq s}.
$$

For arbitrary $p\in \mathbb{F}[\bm{x}]$,  we have the following facts according to equality (\ref{mianformulas1}) and (\ref{2mianformulas2}).
\begin{enumerate}
  \item For fixed $1\leq k\leq \mu$ and $\bm{\alpha} \in
\mathfrak{\Delta}^{{(k)}}$, $(D_{\bm{\rho}^{(k)}}^{\bm{\alpha}}p)({\bm{\xi}^{(k)}})$ can be
linearly expressed by $$\{p({\bm{\xi}^{(k)}+ h\sum_{i=1}^d  \beta_i\bm{\rho}_i^{(k)}}): \bm{\beta}\in \mathfrak{\Delta}^{(k)}\}\cup \{O(h)\}$$ since $\mathfrak{\Delta}^{{(k)}}$ is lower, and moreover, the linear combination
coefficient of each $p({\bm{\xi}^{(k)}+ h\sum_{i=1}^d  \beta_i\bm{\rho}_i^{(k)}})$ is independent of $p\in \mathbb{F}[\bm{x}]$.
  \item For fixed $1\leq l\leq \nu$ and $0\leq m\leq a_1^{(l)}$, $(q_{n(l), m}^{(l)}(D)p)(\bm{\xi}^{(\mu+l)})$ can be
linearly expressed by $$\{p(\bm{\xi}^{(\mu+l)}+\bm{\phi}^{(l)}(r h)): 0\leq r\leq m\}\cup \{O(h)\}.$$ Also, the linear combination
coefficient of each
$p(\bm{\xi}^{(\mu+l)}+\bm{\phi}^{(l)}(r h))$ is independent of $p\in \mathbb{F}[\bm{x}]$.
\end{enumerate}

In brief, we can conclude that there exists a
nonsingular matrix $T$ such that
\begin{equation}\label{juzhen}
\left[\widehat{\bm{\lambda}_h^T \bm{q}}| \widehat{\bm{\lambda}_h^T
f}\right]:=T\left[\bm{\lambda}^T_h \bm{q}| \bm{\lambda}^T_h
f\right]=\left[\bm{\lambda}^T \bm{q}| \bm{\lambda}^T f\right]+\left[E_h|
\bm{\epsilon}_{h}\right],
\end{equation}
where each entry of  $[E_h|\bm{\epsilon}_{h}]$ has the same order as $h$. As a consequence, the linear systems
$$\left(\widehat{\bm{\lambda}_h^T \bm{q}}\right) \bm{x}=\widehat{\bm{\lambda}_h^T
f}\quad \mbox{and} \quad\left(\bm{\lambda}^T_h \bm{q}\right) \bm{x}=\bm{\lambda}^T_h f$$
are equivalent, namely they have the same set of solutions.

(i) From (\ref{juzhen}), it follows that each entry of
matrix $\widehat{\bm{\lambda}_h^T \bm{q}}$ converges to its
corresponding entry of matrix $\bm{\lambda}^T \bm{q}$  as $h$ tends
to zero, which implies that
$$\lim\limits_{h\rightarrow
0}\det\left(\widehat{\bm{\lambda}_h^T \bm{q}}\right)=\det \left(\bm{\lambda}^T
\bm{q}\right).$$
Since $\det (\bm{\lambda}^T \bm{q})\neq 0$, there exists $\eta>0$ such that
$$\det\left(\widehat{\bm{\lambda}_h^T \bm{q}}\right) \neq 0, \quad0<|h|<\eta.$$
Notice that (\ref{juzhen}) directly leads to  $\mathrm{rank}
\left(\widehat{\bm{\lambda}_h^T \bm{q}}\right)=\mathrm {rank}\left(\bm{\lambda}_h^T
\bm{q}\right)$,
$$
\mathrm{ran}P_h=\mathrm{span}_\mathbb{F}\bm{q}, \quad 0<|h|<\eta,
$$
follows, i.e., $\bm{q}$ forms an $\mathbb{F}$-basis for $\mathrm{ran}P_h$. Since $\bm{q}$ is also an $\mathbb{F}$-basis for $\mathrm{ran}P$, we have
$$
\mathrm{ran}P=\mathrm{ran}P_h,\quad 0<|h|<\eta.
$$

(ii) Suppose that $\widetilde{\bm{x}}$ and $\bm{x_0}$ be the unique solutions of nonsingular linear systems
\begin{equation}\label{Leq}
(\bm{\lambda}_h^T \bm{q})\bm{x}=\bm{\lambda}_h^T f
\end{equation}
and
\begin{equation}\label{Heq}
(\bm{\lambda}^T \bm{q})\bm{x}=\bm{\lambda}^T f
\end{equation}
respectively, where $f\in \mathbb{F}[\bm{x}]$ and $0<|h|<\eta$. It is easy to see that
$$
P_h f=\bm{q\widetilde{\bm{x}}}\quad \mbox{and}\quad Pf=\bm{qx_0}.
$$
Notice that, as $h\rightarrow 0$,  $P$ is the pointwise limit of $P_h$ if and only if $P f$ is the coefficientwise limit of $P_h f$ for all $f\in \mathbb{F}[\bm{x}]$. Therefore, it is sufficient to show that for every $f\in \mathbb{F}[\bm{x}]$, the solution vector of system (\ref{Leq}) converges to the one of system (\ref{Heq}) as $h$ tends to zero, namely
$$
\lim_{h\rightarrow 0} \widetilde{\bm{x}}= \bm{x_0}.
$$
By (\ref{juzhen}), the linear system
\begin{equation}\label{widehat}
    \left(\widehat{\bm{\lambda}_h^T
\bm{q}}\right)\bm{x}=\widehat{\bm{\lambda}_h^T f}
\end{equation}
can be rewritten as
$$\left(\bm{\lambda}^T \bm{q}+E_h\right)\bm{x}=\left(\bm{\lambda}^T
f+\bm{\epsilon}_{h}\right).
$$
Since system (\ref{widehat}) is equivalent to system (\ref{Leq}), $\widetilde{\bm{x}}$ is also the unique solution of it. Consequently, using the perturbation analysis of the sensitivity
of linear systems (see for example \cite[p. 80]{Matrixcompu1996}), we have
$$\left\|\widetilde{\bm{x}}-\bm{x_0}\right\|\leq\left\|{(\bm{\lambda}^T \bm{q})}^{-1}\right\|\left\|\bm{\epsilon}_h-E_h\bm{x_0}\right\|+O(h^2).$$
Since each component of vector $\bm{\epsilon}_h-E_h\bm{x_0}$  has the same order as $h$, it follows that $\lim\limits_{h\rightarrow
0}\|\widetilde{ \bm{x}}-\bm{x_0}\|= 0$, or, equivalently,
$\lim\limits_{h\rightarrow 0}\widetilde{\bm{x}}=\bm{x_0}$, which completes the proof of the theorem.
\qed
\end{pf}

Finally, we will give a complete example to
illustrate the conclusions
 of
Theorem \ref{mianthm}.

\begin{exmp}\label{ex1}
Let ${\bm{\xi}^{(1)}}=(1, 1, 1)$, ${\bm{\xi}^{(2)}}=(0, 0, 0)$. Let
$$\bm{q}=(1, x_3, x_2, x_1, x_3^2, x_3 x_2, x_3 x_1)$$
and
\begin{align*}
\bm{\lambda}=\big(&\delta_{\bm{\xi}^{(1)}}, \delta_{\bm{\xi}^{(1)}}
\frac{\partial}{\partial x_1}, \delta_{\bm{\xi}^{(1)}}
\frac{\partial}{\partial x_2}, \delta_{\bm{\xi}^{(1)}}
\frac{\partial}{\partial x_3},
\delta_{{\bm{\xi}^{(2)}}},
\delta_{{\bm{\xi}^{(2)}}}(\frac{\partial }
{\partial x_1}+ \frac{\partial f}{\partial x_2}),\\
& \delta_{{\bm{\xi}^{(2)}}}( \frac{1}{2}\frac{\partial^2 }
{\partial x_1^2}+\frac{\partial^2 }
{\partial x_1 \partial x_2 } +\frac{1}{2}\frac{\partial^2 }
{\partial x_2^2}+\frac{\partial }
{\partial x_2}+ \frac{\partial }{\partial x_3}
)\big)
\end{align*}
be
the $\mathbb{F}$-basis for $\mathrm{ran} P$ and
$\mathrm{ran} P'$, respectively.

From Proposition \ref{D1prop1} and Example 3, we know that
this example is the case of Theorem  \ref{mianthm},
Therefore, we set
$$
\bm{\lambda}_h=(\delta_{(1, 1, 1)}, \delta_{(1+h, 1, 1)}, \delta_{(1, 1+h, 1)}, \delta_{(1, 1, 1+h)},
\delta_{(0, 0, 0)}, \delta_{(h, h^2+h, h^2)}, \delta_{(2h, 4h^2+2h, 4h^2)}).
$$
Recalling the proof of Theorem \ref{mianthm}, we can obtain
$$\widehat{\bm{\lambda}_h^T\bm{q}}=\left(
                                     \begin{array}{ccccccc}
                                       1 & 1 & 1 & 1 & 1 & 1 & 1 \\
                                       0 & 0 & 0 & 1 & 0 & 0 & 1 \\
                                       0 & 0 & 1 & 0 & 0 & 1 & 0 \\
                                       0 & 1 & 0 & 0 & 2+h & 1 & 1 \\
                                       1 & 0 & 0 & 0 & 0 & 0 & 0 \\
                                       0 & h & 1+h & 1 & h^3 & h^2(1+h) & h^2 \\
                                       0 & 1 & 1 & 0 & 7 h^2 & h(7h+3) & 3 h \\
                                      \end{array}
                                   \right).
$$
Then
$$
\det (\widehat{\bm{\lambda}_h^T\bm{q}})=(h-1)(2h+1)(2 h-1)(1+h)^2,
$$
hence
$$
\det (\widehat{\bm{\lambda}_h^T\bm{q}})\neq 0,\quad  0<|h|< \frac{1}{2}.
$$
Consequently, $P$ is the pointwise limit of Lagrange projector $P_h,
0<|h|<\frac{1}{2}$,  as $h$ tends to zero, with the property that
$\mathrm{ran} P'_h=\mathrm{span}_{\mathbb{F}}\bm{\lambda}_h$.

More precisely, we select a test function
$$f(x_1,x_2,x_3)=1+(1-x_1)^2+(1-x_2)^2+(1-x_3)^2$$
to describe the perturbation procedure for the projector in this example.

When $h=1/10, 1/100, 1/1000, \ldots$, we have
\begin{align*}
P_{\frac{1}{10}}f=&4- \frac{34949}{14520}x_3-\frac{439}{7260}x_2-\frac{37867}{9680} x_1-\frac{2303}{2904}x_3^2+\frac{233}{1452}x_3 x_2+\frac{7767}{1936}x_3 x_1,\\
P_{\frac{1}{100}}f=&4- \frac{2600449499}{1274614950}x_3-\frac{46747801}{2549229900}x_2-\frac{483294631}{121391900} x_1-\frac{24977753}{25492299}x_3^2\\
&+\frac{722401}{25492299}x_3 x_2+\frac{9690171}{2427838}x_3 x_1,\\
P_{\frac{1}{1000}}f=&4- \frac{251000494994999}{125249623999500}x_3-\frac{496749753001}{250499247999000}x_2-\frac{333833081249167}{83499749333000} x_1\\
&-\frac{249997752503}{250499247999}x_3^2
+\frac{747249001}{250499247999}x_3 x_2+\frac{667833161997}{166999498666}x_3 x_1,\\
\cdots&\\
P f=&4-2 x_3-4 x_1-x_3^2+4 x_3 x_1.
\end{align*}
\qed
\end{exmp}

\bibliographystyle{elsarticle-num}
\bibliography{ref}

\end{document}